\documentclass[reqno]{article}
\usepackage{amssymb,amsmath}
\usepackage{amsbsy}
\usepackage{color}
\usepackage[latin1]{inputenc}

\newcommand{\disp}{\displaystyle}

\def\Og{\mathcal{O}}
\def\Om{{\Omega}}

\def\Ombar{\overline{\Om}}
\def\jnt{\disp\int}
\def\jjnt{\jnt\!\!\!\!\jnt}

\def\ld{{\lambda}}

\def\al{{\alpha}}

\def\R{{\mathbb{R}}}

\def\eps{\varepsilon}
\newcommand{\Fin}{\hfill$\Box$}

\numberwithin{equation}{section}
\newtheorem{lemma}{Lemma}[section]
\newtheorem{proposition}{Proposition}[section]
\newtheorem{theorem}{Theorem}[section]
\newtheorem{corollary}{Corollary}[section]

\begin{document}

\title{Null controllability for a parabolic-elliptic coupled system}

\author
{
E. Fern\'andez-Cara\thanks{Dpto.~E.D.A.N., Universidad de Sevilla, Aptdo.~1160, 41080~Sevilla, Spain, {\tt cara@us.es}.},
\quad J. Limaco\thanks{Inst.~Matem\'atica, Universidade Federal Fluminense, Valonguinho, 24020-140,  Niter\'oi, RJ, Brasil, {\tt jlimaco@vm.uff.br}.} ,
\quad S. B. de Menezes\thanks{Corresponding author: Dpto.~Matem\'atica, Universidade Federal do Cear\'a, Campus do Pici - Bloco 914, 60455-760, Fortaleza, CE, Brasil, {\tt silvano@mat.ufc.br}.}
}


\maketitle

\begin{abstract}
   In this paper, we prove the null controllability of some parabolic-elliptic systems.
   The control is distributed, locally supported in space and appears only in one PDE.
   The arguments rely on fixed-point reformulation and suitable Carleman estimates for the solutions to the adjoint system.
   Under appropriate assumptions, we also prove that the solution can be obtained as the asymptotic limit of some similar parabolic systems.

\end{abstract}

\

\noindent{\bf Mathematics Subject Classification 2000:} 35B37, 35A05, 35B40.

\noindent{\bf Keywords:} Null controllability, parabolic-elliptic systems, Carle\-man inequalities.

\section{Introduction and main results}\label{introd}


   Let $\Om$ be a bounded domain of $\R^{N}$ ($N\geq1$), with boundary $\Gamma=\partial \Om$ of class $C^2$.
   We fix $T>0$ and we denote by $Q$ the cylinder $Q=\Om\times (0,T)$, with lateral boundary $\Sigma=\Gamma\times (0,T)$.
   We also consider a non-empty
   (small) open set $\Og \subset \Om$;
   as usual, $1_\Og$ denotes the characteristic function of $\Og$.

   Throughout this paper, $C$ (and sometimes $C_{0}$, $K$, $K_0$, \dots) denotes various positive constants.
   Frequently, we will emphasize the fact that $C$ depends on
   (say) $f$ by writing $C = C(f)$.
   The inner product and norm in~$L^{2}(\Om)$ will be denoted, respectively, by $(\cdot\,,\cdot)$ and~$\|\cdot\|$.
   On the other hand, $\|\cdot\|_\infty$ will stand for the norm in~$L^\infty(Q)$.

   We will consider the following semilinear parabolic-elliptic coupled systems
   \begin{equation}\label{1.1}
\left\{
\begin{array}{l}
\disp y_{t} - \Delta y = F(y,z)+ v1_{\Og}\;\;\mbox{in}\;\; Q,\\
\disp -\Delta z= f(y,z)\;\;\mbox{in}\;\; Q,\\
\disp y(x,t)=0,\; z(x,t)=0 \;\;\mbox{on}\;\; \Sigma,\\
\disp y(x,0)=y^{0}(x)\;\;\mbox{in}\;\; \Om
\end{array}
\right.
   \end{equation}
and
   \begin{equation}\label{1.1p}
\left\{
\begin{array}{l}
\disp y_{t} - \Delta y = F(y,z) \;\;\mbox{in}\;\; Q,\\
\disp -\Delta z= f(y,z) + w1_{\Og} \;\;\mbox{in}\;\; Q,\\
\disp y(x,t)=0,\; z(x,t)=0 \;\;\mbox{on}\;\; \Sigma,\\
\disp y(x,0)=y^{0}(x)\;\;\mbox{in}\;\; \Om,
\end{array}
\right.
   \end{equation}
under some hypotheses for $F$ and $f$.

   In~(\ref{1.1}) and~(\ref{1.1p}), we have $y=y(x,t)$ and~$z=z(x,t)$;
   $1_{\Og}$ is the characteristic function of $\Og$ and $y^0 = y^0(x)$ is the initial state.
   We will assume that the possibly nonlinear functions $F:\R\times \R \mapsto \R$ and $ f:\R \mapsto \R$ satisfy:
   \begin{equation}\label{1.2}
\left\{
\begin{array}{l}
\disp \text{$F$ and $f$ are (globally) Lipschitz-continuous,} \\
\noalign{\smallskip}
\disp F(0,0) = f(0,0) = 0, \quad {\partial f \over \partial z} (y,z) \leq \mu < \mu_1 \ \ \text{a.e.},\\
\end{array}
\right.
   \end{equation}
where $\mu_1$ the first eigenvalue of the Dirichlet Laplacian in $\Om$.

   The analysis of systems of the kind \eqref{1.1} and~\eqref{1.1p} can be justified by several applications.
   Let us indicate two of them:

\begin{itemize}

\item Reaction-diffusion systems with origin in physics, chemistry, biology, etc.\ where two scalar ``populations'' interact and the natural time scale of the growth rate is much smaller for one of them than for the other one.
   Precise examples can be found in the study of prey-predator interaction, chemical heating, tumor growth therapy, etc.

\item Semiconductor modeling, where one of the state variables is
   (for example) the density of holes and the other one is the electrical potential of the device;
   see for instance~\cite{Jungel}.
   Other problems with this motivation will be analyzed with more detail by the authors in the next future.

\end{itemize}

   The system (\ref{1.1})
   (resp.~\eqref{1.1p}) is well-posed in the sense that, for each $y^0\in L^2(\Om)$ and each $v \in L^2(\Og\times(0,T))$
   (resp.~$w \in L^2(\Og\times(0,T))$) possesses exactly one solution $(y,z)$, with
   $$
 y \in L^2(0,T;H_0^1(\Om)), \quad y_t \in L^2(0,T;H^{-1}(\Om)), \quad z \in L^2(0,T;D(-\Delta)).
   $$
   This statement is justified in~Appendix~A.

   In this paper we will analyze some controllability properties of~\eqref{1.1} and~\eqref{1.1p}.

   It will be said that (\ref{1.1})
   (resp.~\eqref{1.1p}) is {\it null-controllable} at time $T$ if the following holds:
   for any given $y^0 \in L^2(\Om)$, there exist controls $v \in L^2(\Og\times(0,T))$
   (resp.~controls $w \in L^2(\Og\times(0,T))$) and associated solutions satisfying
   $$
z \in C^0([0,T];L^2(\Om))
   $$
and
   \begin{equation}\label{eq2n}
 y(x,T)=0 \ \mbox{ in } \ \Om, \quad \limsup_{t \to T^-} \|z(\cdot\,,t)\| = 0,
   \end{equation}
with an estimate of the form
   \begin{equation}
\|v\|_{L^2(\Og\times(0,T))}\leq C \|y^0\|
\quad (\mbox{resp.\ }\|w\|_{L^2(\Og\times(0,T))}\leq C \|y^0\|). \label{eq3n}
   \end{equation}

   This inequality indicates that the ``null controls'' can be chosen depending continuously on the initial data.

   The control of PDEs equations and systems has been the subject of a lot of papers the last years.
   In particular, important progress has been made recently in the controllability analysis of semi-linear parabolic equations.
   We refer to the works \cite{Coron-book, D-C-B-Z, FPZ, FCZ2, FI, Oleg1, Z0, EZ-rev} and the references therein.
   Consequently, it is natural to try to extend the known results to systems of the kind~\eqref{1.1} and~\eqref{1.1p}.


   The main results in this paper are the following:

\begin{theorem} \label{t1n}
   Assume that
   \begin{equation}
\text{$F(y,z) \equiv F_0(y) + bz$,  with $F_0$ Lipschitz-continuous, $b \in \R$,} \label{1.2a}
   \end{equation}
   \begin{equation}
\text{$f(y,z) \equiv cy + dz$,  with $c,d \in \R$, $c \not=0$, $d < \mu_1$.} \label{1.2b}
   \end{equation}
   Then \eqref{1.1} is null-controllable at any time~$T>0$.
\end{theorem}

\begin{theorem} \label{t2n}
   Let us assume that \eqref{1.2a} holds and
    \begin{equation}
\text{$f(y,z) \equiv f_0(y) + dz$,  with $f_0$ Lipschitz-continuous, $d \in \R$, $d < \mu_1$.} \label{1.2c}
   \end{equation}
   Then \eqref{1.1p} is null-controllable at any time~$T>0$.
\end{theorem}

   The proofs of these results rely on relatively well known arguments and some new estimates.

   More precisely, in a first step, we will first consider similar {\it linearized} systems of the form
   \begin{equation}\label{1.1-lin}
\left\{
\begin{array}{l}
\disp y_{t} - \Delta y = a(x,t)y + bz + v1_{\Og}\;\;\mbox{in}\;\; Q,\\
\disp -\Delta z= cy + dz \;\;\mbox{in}\;\; Q,\\
\disp
 y(x,t)=0,\; z(x,t)=0 \;\;\mbox{on}\;\; \Sigma,\\
\disp
 y(x,0)=y^{0}(x)\;\;\mbox{in}\;\; \Om
\end{array}
\right.
   \end{equation}
and
   \begin{equation}\label{1.1p-lin}
\left\{
\begin{array}{l}
\disp y_{t} - \Delta y = a(x,t)y + bz \;\;\mbox{in}\;\; Q,\\
\disp -\Delta z= cy + dz + w1_{\Og} \;\;\mbox{in}\;\; Q,\\
\disp
 y(x,t)=0,\; z(x,t)=0 \;\;\mbox{on}\;\; \Sigma,\\
\disp
 y(x,0)=y^{0}(x)\;\;\mbox{in}\;\; \Om.
\end{array}
\right.
   \end{equation}
   We will establish null controllability results for~\eqref{1.1-lin} and~\eqref{1.1p-lin} by previously proving appropriate Carleman estimates for the solutions to the associated {\it adjoint} systems.
   Then, in a second step, we will adapt a fixed-point argument to get the null controllability results stated in Theorems~\ref{t1n} and~\ref{t2n}.

\

   In this paper, we will also consider systems of the form
   \begin{equation}\label{1epsn-a}
\left\{
\begin{array}{l}
 \disp y_{t} - \Delta y = F(y,z)\;\;\mbox{in}\;\; Q,\\
\disp \eps z_t -\Delta z= f(y,z)+ w1_{\Og}\;\;\mbox{in}\;\; Q,\\
\disp
 y(x,t)=0,\; z(x,t)=0 \;\;\mbox{on}\;\; \Sigma,\\
\disp
 y(x,0)=y^{0}(x), \,\,z(x,0)=z^0(x)\;\;\mbox{in}\;\; \Omega.
\end{array}
\right.
   \end{equation}

   It will be shown that, under the assumptions of~Theorem~\ref{cnepsl}, \eqref{1epsn-a} is uniformly null-controllable as $\eps \to 0$, i.e.~null-controlable with controls $w_\eps$ satisfying the estimates \eqref{eq3n} with $C$ independent of~$\eps$.
   We will also see that the $w_\eps$ can be chosen in such a way that they converge weakly to a null control of~\eqref{1.1p}
   (see~Theorem~\eqref{TFep} below).

\

   This paper is organized as follows.
   In~Section~\ref{Sec2}, we introduce some adjoint
   (backwards in time) parabolic-elliptic systems and we prove that their solutions satisfy suitable Carleman estimates.
   In~Section~\ref{Sec3}, we deduce from these estimates null controllability results for ~\eqref{1.1-lin} and~\eqref{1.1p-lin}.
   Section~\ref{Sec4} deals with the proofs of Theorems~\ref{t1n} and~\ref{t2n}.
   The uniform null controllability property of~\eqref{1epsn-a} and the convergence of the associated null controls are established in~Section~\ref{Sec5}.
   Finally, we give the proofs of some technical results in~Section~\ref{Sec6}
   (Appendix~A).


\section{Some Carleman estimates}\label{Sec2}

   We will first consider the general linear backwards in time system
   \begin{equation}\label{as-0}
\left\{
\begin{array}{l}
 - \varphi_t-\Delta \varphi=a(x,t)\varphi+c(x,t)\psi\;\;\mbox{in}\;\; Q,\\
\disp -\Delta \psi=b(x,t)\varphi+d(x,t)\psi\;\;\mbox{in}\;\; Q,\\
\disp \varphi=0, \ \psi=0 \;\;\mbox{on}\;\; \Sigma,\\
\disp \varphi(x,T)=\varphi^T(x)\;\;\mbox{in}\;\; \Om,
\end{array}
\right.
   \end{equation}
where $\varphi^T \in L^2(\Om)$ and we assume that
   \begin{equation}\label{hyp-1.1-0}
a,b,c,d \in L^\infty(Q), \ \ d \leq \mu < \mu_1 \ \text{a.e.}
   \end{equation}

   Also, it will be convenient to introduce a new non-empty open set $\Og_0$, with $\Og_0 \Subset \Og$.
   We will need the following result, due to Fursikov and Imanuvilov~\cite{FI}:

\begin{lemma}\label{lem2.1}
   There exists a function $\alpha_0 \in C^2(\overline{\Om})$ satisfying:
   $$
\left\{
\begin{array}{l}
\alpha_0 (x)>0 \quad \forall x \in \Om, \quad \alpha_0=0 \quad \forall  x \in \partial \Om,\\
|\nabla \alpha_0(x)|>0 \quad \forall x \in \Ombar\setminus \Og_0.
\end{array}
\right.
   $$
\end{lemma}

   Let us introduce the functions
   $$
\beta(t) := t(T-t), \ \phi (x,t) := \frac{e^{\lambda\alpha_0(x)}}{\beta(t)}\,, \
\overline\alpha(x) := e^{k\lambda}-e^{\lambda\alpha_0(x)}, \ \alpha(x,t) := \frac{\overline\alpha(x)}{\beta(t)} \,,
   $$
where $k> \|\alpha_0\|_{L^\infty}+\log 2$ and $\lambda >0$.
   Also, let us set
   $$
\begin{array}{l}
\disp \hat\alpha(t) := \min_{x \in \overline\Om} \alpha(x,t), \quad \alpha^*(t) := \max_{x \in \overline\Om} \alpha(x,t), \\
\disp \hat\phi(t) := \min_{x \in \overline\Om} \phi(x,t), \quad \phi^*(t) := \max_{x \in \overline\Om} \alpha(x,t).
\end{array}
   $$

   Then the following Carleman estimates hold:

\begin{proposition}\label{t2-0}
   Assume that \eqref{hyp-1.1-0} holds.
   There exist positive constants $\lambda_0$, $s_0$ and $C_0$ such that, for any $s \geq s_0$ and $\lambda \geq \lambda_0$ and any $\varphi^T \in L^2(\Om)$, the associated solution to~\eqref{as-0} satisfies
   \begin{equation}\label{cp1}
\begin{array}{c}
\displaystyle \jjnt_Q \!\!e^{-2s\alpha}\!\left[(s\phi)^{-1}\!\left(|\varphi_t|^2\!+\!|\Delta \varphi|^2 \right)\!+\! \lambda^2 (s\phi) |\nabla \varphi|^2 \!+\! \lambda^4(s\phi)^3|\varphi|^2\right]dxdt \\
\displaystyle \leq C_0 \left( \jjnt_Q e^{-2s\alpha}|\psi|^2 + \jjnt_{\Og_0\times(0,T)} e^{-2s\alpha} \lambda^4(s\phi)^3|\varphi|^2  \right)dxdt
\end{array}
   \end{equation}
and
   \begin{equation}\label{ce1}
\begin{array}{c}
\displaystyle \jjnt_Q \!\!e^{-2s\alpha}\!\left[(s\phi)^{-1}\!|\Delta \psi|^2\!+\! \lambda^2 (s\phi) |\nabla \psi|^2 \!+\! \lambda^4(s\phi)^3|\psi|^2\right]dxdt\!  \\
\displaystyle \leq C_0 \left( \jjnt_Q e^{-2s\alpha}|\varphi|^2 + \jjnt_{\Og_0\times(0,T)} e^{-2s\alpha} \lambda^4(s\phi)^3|\psi|^2  \right)dxdt
\end{array}
   \end{equation}
Furthermore, $C_0$ and $\lambda_0$ only depend on~$\Om$ and~$\Og$ and $s_0$ can be chosen of the form
   \begin{equation}\label{2.6p}
\displaystyle s_0 = \sigma_0 (T + T^2),
   \end{equation}
where $\sigma_0$ only depends on~$\Om$, $\Og$, $\|a\|_\infty$, $\|b\|_\infty$, $\|c\|_\infty$ and~$\|d\|_\infty$.
\end{proposition}

   This result is proved in~\cite{FI}.
   In fact, similar Carleman inequalities are established there for more general linear parabolic equations.
   The explicit dependence in time of the constants is not given in~\cite{FI}.
   We refer to~\cite{Cara-Guerrero}, where the above formula for $s_0$ is obtained.

   For further purpose, we introduce the following notation:
   $$
\begin{array}{l}
\disp I(s,\lambda; \varphi) = \jjnt_Q \!\!e^{-2s\alpha}\!\left[(s\phi)^{-1}\!\left(|\varphi_t|^2\!+\!|\Delta \varphi|^2 \right)\!+\! \lambda^2 (s\phi) |\nabla \varphi|^2 \!+\! \lambda^4(s\phi)^3|\varphi|^2\right]dxdt
\end{array}
   $$
and
   $$
\begin{array}{l}
\disp \tilde{I}(s,\lambda; \psi) = \jjnt_Q \!\!e^{-2s\alpha} \!\left[ (s\phi)^{-1} |\Delta \psi|^2 \!+\! \lambda^2 (s\phi) |\nabla \psi|^2 \!+\! \lambda^4(s\phi)^3|\psi|^2\right]dxdt .
\end{array}
   $$

   Now, we will deduce several consequences from Proposition~\ref{t2-0} under particular hypotheses on the coefficients of~\eqref{as-0}.
   First, it will be assumed that $c$ is a.e.~equal to a non-zero constant and $b$ and $d$ do not depend of $t$:
   \begin{equation}\label{hyp-1.1}
a \in L^\infty(Q), \ c \in \R, \ c \not= 0,\ b,d \in L^\infty(\Om),\ d \leq \mu < \mu_1 \ \text{a.e.}
   \end{equation}

   Accordingly, \eqref{as-0} reads:
   \begin{equation}\label{as}
\left\{
\begin{array}{l}
 - \varphi_t-\Delta \varphi=a(x,t)\varphi+c\psi\;\;\mbox{in}\;\; Q,\\
\disp -\Delta \psi=b(x)\varphi+d(x)\psi\;\;\mbox{in}\;\; Q,\\
\disp \varphi(x,t) = 0, \ \psi(x,t) = 0 \;\;\mbox{on}\;\; \Sigma,\\
\disp \varphi(x,T)=\varphi^T(x)\;\;\mbox{in}\;\; \Om.
\end{array}
\right.
   \end{equation}

\begin{proposition}\label{t2}
   Assume that \eqref{hyp-1.1} holds.
   There exist positive constants $\lambda_0$, $s_0$ and $C_1$ such that, for any $s \geq s_0$ and $\lambda \geq \lambda_0$ and any $\varphi^T \in L^2(\Om)$, the associated solution to~\eqref{as} satisfies
   \begin{equation}\label{Carleman-1}
I(s,\lambda; \varphi) + \tilde{I}(s,\lambda; \psi) \leq C_1 \jjnt_{\Og\times(0,T)} e^{-4s\hat\alpha + 2s\alpha^*} \lambda^8 (s\phi^*)^7 |\varphi|^2dxdt .
   \end{equation}
   Furthermore, $C_1$ and $\lambda_0$ only depend on~$\Om$ and~$\Og$ and $s_0$ can be chosen of the form
   \begin{equation}\label{2.6p}
\displaystyle s_1 = \sigma_1 (T + T^2),
   \end{equation}
where $\sigma_1$ only depends on~$\Om$, $\Og$, $\|a\|_\infty$, $\|b\|_{L^\infty}$, $|c|$ and~$\|d\|_{L^\infty}$.
\end{proposition}

\noindent
{\sc Proof:}
   Obviously, it will be sufficient to show that there exist $\lambda_0$, $s_0$ and $C_1$ such that, for any small $\eps > 0$, one has:
   \begin{equation}\label{4l}
\begin{array}{l}
\disp I(s,\lambda; \varphi) + \tilde{I}(s,\lambda; \psi) \leq C\eps I(s,\lambda; \varphi) + C\eps \tilde{I}(s,\lambda; \psi)  \\
\disp \qquad + \ C_1 \jjnt_{\Og\times(0,T)} e^{-4s\hat\alpha + 2s\alpha^*} \lambda^8 (s\phi^*)^7 |\varphi|^2dxdt .
\end{array}
   \end{equation}

   We start from \eqref{cp1} and~\eqref{ce1}.
   After addition, by taking $\sigma_1$ sufficiently large and $s \geq \sigma_1 (T + T^2)$, we obtain:
   \begin{equation}\label{5l}
\begin{array}{l}
\disp I(s,\lambda; \varphi) + \tilde{I}(s,\lambda; \psi) \\
\disp \quad \leq C \jjnt_{\Og_0\times(0,T)} e^{-2s\alpha} \lambda^4(s\phi)^3\left(|\varphi|^2+|\psi|^2\right)dxdt.
\end{array}
   \end{equation}

   Let us now introduce a function $\xi \in \mathcal{D}(\Og)$ satisfying $0<\xi \leq 1$ and $\xi \equiv 1$ in~$\Og_0$.
   Then
   \begin{equation}
\begin{array}{l}
\disp \jjnt_{\Og_0\times(0,T)} e^{-2s\alpha} \lambda^4(s\phi)^3|\psi|^2dxdt \leq \jjnt_{\Og\times(0,T)} e^{-2s\alpha} \lambda^4(s\phi)^3 \xi |\psi|^2 dxdt\\
\disp \quad = \jjnt_{\Og\times(0,T)} e^{-2s\alpha} \lambda^4(s\phi)^3 \xi(x) \psi\left(-\frac{1}{c}(\varphi_t+\Delta \varphi+a(x,t)\varphi)\right)dxdt \\
\disp \quad = - \jjnt_{\Og\times(0,T)} e^{-2s\alpha} \lambda^4(s\phi)^3 \frac{\xi(x)}{c} \,\psi \,\varphi_t \;dxdt\\
\disp \phantom{\quad = } - \jjnt_{\Og\times(0,T)} e^{-2s\alpha} \lambda^4(s\phi)^3 \frac{\xi(x)}{c} \,\psi \,\Delta\varphi\; dxdt\\
\disp \phantom{\quad = } - \jjnt_{\Og\times(0,T)} e^{-2s\alpha} \lambda^4(s\phi)^3 \frac{\xi(x)}{c} \,a(x,t) \,\psi \,\varphi\; dxdt \\
\disp \quad := M_1+M_2+M_3.
\label{6l}
\end{array}
   \end{equation}

   Let us compute and estimate the $M_i$.

   First,
   $$
\begin{array}{l}
\disp M_1 = \jjnt_{\Og\times(0,T)} e^{-2s\alpha} \frac{2\xi(x)}{c} \lambda^4 s^4 \phi^3 \alpha_t \psi \varphi\;dxdt \\
\disp \phantom{M_1 = } + \jjnt_{\Og\times(0,T)} e^{-2s\alpha} \frac{3\xi(x)}{c} \lambda^4 s^3 \phi^2 \phi_t \psi \varphi\;dxdt \\
\disp \phantom{M_1 = } + \jjnt_{\Og\times(0,T)} e^{-2s\alpha} \frac{\xi(x)}{c} \lambda^4 (s\phi)^3 \psi_t \varphi \;dxdt.
\end{array}
   $$
   Using that $|\alpha_t| \leq C\phi^2$ and $|\phi_t| \leq C\phi^2$ for some $C>0$, we get:
   $$
\begin{array}{l}
\disp M_1 \leq C \jjnt_{\Og\times(0,T)} e^{-2s\alpha} \lambda^4 s^4 \phi^5 |\psi| \, |\varphi|\;dxdt \\
\disp \phantom{M_1 \leq } + \jjnt_{\Og\times(0,T)} e^{-2s\alpha} \lambda^4 (s\phi)^3 |\psi_t| \, |\varphi|\;dxdt \\
\disp \phantom{M_1 } \leq \eps \tilde{I}(s,\lambda; \psi) + C_\eps \jjnt_{\Og\times(0,T)} e^{-2s\alpha} \lambda^4 s^5 \phi^7 |\varphi|^2 dxdt \\
\disp \phantom{M_1 \leq } + \jjnt_{\Og\times(0,T)} e^{-2s\alpha} \lambda^4 (s\phi)^3 |\psi_t| \, |\varphi|\;dxdt.
\end{array}
   $$
   The last integral in this inequality can be bounded as follows:
   $$
\begin{array}{l}
\disp \jjnt_{\Og\times(0,T)} e^{-2s\alpha} \lambda^4 (s\phi)^3 |\psi_t| \, |\varphi|\;dxdt \\
\disp \quad \leq \jjnt_{\Og\times(0,T)} e^{-2s\hat\alpha} \lambda^4 (s\phi^*)^3 |\psi_t| \, |\varphi|\;dxdt \\
\disp \quad = \jnt_0^T e^{-2s\hat\alpha(t)} \lambda^4 (s\phi^*(t))^3 \|\psi_t(\cdot\,,t)\|_{L^2(\Og)}\|\varphi(\cdot\,,t)\|_{L^2(\Og)} \;dt\\
\disp \quad \leq C \jnt_0^T e^{-2s\hat\alpha(t)} \lambda^4 (s\phi^*(t))^3 \|\varphi_t(\cdot\,,t)\| \|\varphi(\cdot\,,t)\|_{L^2(\Og)}\;dt \\
\disp \quad = C \jnt_0^T e^{-s\alpha^*} (s\phi^*(t))^{-1/2} \|\varphi_t(\cdot\,,t)\| \cdot e^{-2s\hat\alpha+s\alpha^*} \lambda^4 (s\phi^*)^{7/2} \|\varphi(\cdot\,,t)\|_{L^2(\Og)} \;dt\\
\disp \quad \leq \eps I(s,\lambda; \varphi) + C_\eps \jjnt_{\Og\times(0,T)} e^{-4s\hat\alpha+2s\alpha^*} \lambda^8 (s\phi^*)^7 |\varphi|^2 dxdt.
\end{array}
   $$
   Here, we have used that
   $$
\left\{
\begin{array}{l}
\disp -\Delta\psi_t = b\varphi_t + d \psi_t \ \text{ in }\ \Om, \\
\disp \psi_t = 0 \ \text{ on } \ \partial\Om,
\end{array}
\right.
   $$
whence we obviously need $b$ and $d$ independent of $t$.

   Thus, the following is found:
   \begin{equation}\label{estim-M1}
\begin{array}{l}
\disp M_1 \leq \eps I(s,\lambda; \varphi) + \eps \tilde{I}(s,\lambda; \psi) \\
\disp \phantom{M_1} + C_\eps \jjnt_{\Og\times(0,T)} e^{-4s\hat\alpha+2s\alpha^*} \lambda^8 (s\phi^*)^7 |\varphi|^2dxdt.
\end{array}
   \end{equation}

   Secondly, we see that
   \begin{equation}\label{estim-M2}
\begin{array}{l}
\disp M_2 = -\jjnt_{\Og\times(0,T)} \Delta\left( e^{-2s\alpha}  \lambda^4 (s\phi)^3 \frac{\xi(x)}{c} \psi \right) \varphi\;dxdt \\
\disp \phantom{M_2 \leq } C \jjnt_{\Og\times(0,T)} e^{-2s\alpha}
\left[ \lambda^6 (s\phi)^5 |\psi| + \lambda^5 (s\phi)^4 |\nabla\psi| + \lambda^4 (s\phi)^3 |\Delta\psi| \right] \varphi\;dxdt \\
\disp \phantom{M_2} \leq \eps \tilde{I}(s,\lambda; \psi) + C_\eps \jjnt_{\Og\times(0,T)} e^{-2s\alpha} \lambda^8 (s\phi)^7 |\varphi|^2 dxdt.
\end{array}
   \end{equation}
   Here, we have used the identity
   $$
\begin{array}{c}
\disp \Delta\left( e^{-2s\alpha} \phi^3 \frac{\xi(x)}{c} \psi \right) = \Delta\left( e^{-2s\alpha} \phi^3 \frac{\xi(x)}{c} \right) \psi \\
\disp + 2 \nabla\left( e^{-2s\alpha} \phi^3 \frac{\xi(x)}{c} \right) \cdot \nabla\psi + e^{-2s\alpha} \phi^3 \frac{\xi(x)}{c} \Delta\psi
\end{array}
   $$
and the estimates
   $$
|\Delta(e^{-2s\alpha} \phi^3 \frac{\xi(x)}{c})| \leq C e^{-2s\alpha} \lambda^2 s^2 \phi^5
\ \text{ and } \
|\nabla(e^{-2s\alpha} \phi^3 \frac{\xi(x)}{c})| \leq C e^{-2s\alpha} \lambda s \phi^4.
   $$

   Finally, it is immediate that
   \begin{equation}\label{estim-M3}
M_3 \leq \eps \tilde{I}(s,\lambda; \psi) + C_\eps \jjnt_{\Og\times(0,T)} e^{-2s\alpha} \lambda^4 (s\phi)^3 |\varphi|^2 dxdt.
   \end{equation}

   From \eqref{5l}, \eqref{6l} and~\eqref{estim-M1}--\eqref{estim-M3}, we directly obtain \eqref{4l} for all small $\eps>0$.
   This ends the proof.
\Fin

\

   An almost immediate consequence of Proposition~\ref{t2} is the following observability inequality:

\begin{corollary}\label{corol2}
   Under the assumptions of Proposition~\ref{t2}, there exists constants $M$ and~$K$, depending on $\Om$, $\Og$, $T$, $\|a\|_\infty$, $\|b\|_{L^\infty}$, $|c|$ and~$\|d\|_{L^\infty}$, such that every solution $(\varphi,\psi)$ to \eqref{as} verifies:
   \begin{equation}\label{do}
\begin{array}{l}
\displaystyle \|\varphi(\cdot\,,0)\|^2 + \|\psi(\cdot\,,0)\|^2 \leq M \jjnt_{\Og\times(0,T)}e^{-2K/(T-t)} |\varphi|^2dxdt.
\end{array}
   \end{equation}
\end{corollary}

\noindent
{\sc Proof:}
   From the Carleman inequality in~Proposition~\ref{t2} with $s = s_1$ and~$\lambda = \lambda_1$, we see that there exists a constant $C>0$ such that
   $$
\jjnt_Q e^{-2s\alpha} \phi^3 \left(|\varphi|^2+|\psi|^2\right) dxdt \leq C \jjnt_{\Og\times(0,T)} (t(T-t))^{-7} e^{-4s\hat\alpha+2s\alpha^*} |\varphi|^2 dxdt.
   $$
   Since $2\hat\alpha-\alpha^* \equiv a^*/\beta(t)$ for some $a^* > 0$ and $e^{-2s\alpha} \phi^3$ is uniformly bounded from below in $\Om \times [T/4,3T/4]$ we find that
   \begin{equation}\label{20l}
\jnt_{T/4}^{3T/4}\int_\Om \left(|\varphi|^2+|\psi|^2\right) dxdt \leq C \jjnt_{\Og\times(0,T)} e^{-2K/(T-t)} |\varphi|^2 dxdt
   \end{equation}
for some $C, K > 0$.

   On the other hand, we can easily get from \eqref{as} the standard
   (backwards) energy inequalities
   \begin{equation}\label{21l}
-\frac{1}{2} \frac{d}{dt} \|\varphi\|^2 + \|\varphi\|_{H_0^1}^2 \leq C \left (\|\varphi\|^2 + \|\psi\|^2 \right) \ \text{ and } \ \|\psi\|_{H_0^1}^2 \leq C \|\varphi\|^2.
   \end{equation}
   This yields
   \begin{equation}
-\frac{1}{2} \frac{d}{dt} \|\varphi\|^2 \leq C \|\varphi\|^2,
   \end{equation}
whence we deduce that
   \begin{equation}\label{24l}
\|\varphi(\cdot\,,0)\|^2 \leq C \|\varphi(\cdot\,,t)\|^2 \;\;\mbox{ for all } t.
   \end{equation}

   Combining \eqref{20l}, \eqref{24l} and the second part of \eqref{21l}, we obtain at once \eqref{do}.
\Fin

\

   Now, we will assume that $b$ is a non-zero constant:
   \begin{equation}\label{hyp-1.1p}
a,c,d \in L^\infty(Q), \ b \in \R, \ b \not= 0,\ d \leq \mu < \mu_1 \ \text{a.e.}
   \end{equation}

   The corresponding \eqref{as-0} becomes
   \begin{equation}\label{asp}
\left\{
\begin{array}{l}
 - \varphi_t-\Delta \varphi=a(x,t)\varphi+c(x,t)\psi\;\;\mbox{in}\;\; Q,\\
\disp -\Delta \psi=b\varphi+d(x,t)\psi\;\;\mbox{in}\;\; Q,\\
\disp \varphi(x,t) = 0, \ \psi(x,t) = 0 \;\;\mbox{on}\;\; \Sigma,\\
\disp \varphi(x,T)=\varphi^T(x)\;\;\mbox{in}\;\; \Om.
\end{array}
\right.
   \end{equation}

\begin{proposition}\label{t2p}
   Assume that \eqref{hyp-1.1p} holds.
   There exist positive constants $\lambda_2$, $s_2$ and $C_2$ such that, for any $s \geq s_2$ and $\lambda \geq \lambda_2$ and any $\varphi^T \in L^2(\Om)$, the associated solution to~\eqref{asp} satisfies
   $$
I(s,\lambda; \varphi) + \tilde{I}(s,\lambda; \psi) \leq C_2 \jjnt_{\Og\times(0,T)} e^{-2s\alpha} \lambda^8 (s\phi)^7 |\psi|^2 dxdt.
   $$
   Furthermore, $C_2$ and $\lambda_2$ only depend on~$\Om$ and~$\Og$ and $s_2$ can be chosen of the form
   $$
\displaystyle s_2 = \sigma_2 (T + T^2),
   $$
where $\sigma_2$ only depends on~$\Om$, $\Og$, $\|a\|_\infty$, $|b|$, $\|c\|_\infty$ and~$\|d\|_\infty$.
\end{proposition}

\noindent
{\sc Proof:}
   We start again from \eqref{5l}.
   Recalling that $\xi \in \mathcal{D}(\Og)$, $0<\xi \leq 1$ and $\xi \equiv 1$ in~$\Og_0$, we see that
   \begin{equation}
\begin{array}{l}\label{new-p}
\disp \jjnt_{\Og_0\times(0,T)} e^{-2s\alpha} \lambda^4(s\phi)^3|\varphi|^2 dxdt\leq \jjnt_{\Og\times(0,T)} e^{-2s\alpha} \lambda^4(s\phi)^3 \xi |\varphi|^2 dxdt \\
\disp \quad = \jjnt_{\Og\times(0,T)} e^{-2s\alpha} \lambda^4(s\phi)^3 \xi(x) \varphi\left(-\frac{1}{b}(\Delta \psi+d(x,t)\psi)\right)dxdt \\
\disp \quad = - \jjnt_{\Og\times(0,T)} e^{-2s\alpha} \lambda^4(s\phi)^3 \frac{\xi(x)}{b} \varphi \Delta\psi\;dxdt \\
\disp \phantom{\quad = } - \jjnt_{\Og\times(0,T)} e^{-2s\alpha} \lambda^4(s\phi)^3 \frac{\xi(x)}{b} d(x,t) \varphi \psi \;dxdt\\
\disp \quad := M'_1+M'_2.
\end{array}
   \end{equation}

   As in the proof of Proposition~\ref{t2}, it is not difficult compute and estimate the $M'_i$.
   Indeed,
   \begin{equation}\label{estim-M1p}
\begin{array}{l}
\disp M'_1 = -\jjnt_{\Og\times(0,T)} \Delta\left( e^{-2s\alpha}  \lambda^4 (s\phi)^3 \frac{\xi(x)}{b} \varphi \right) \psi \;dxdt\\
\disp \phantom{M'_1}  \leq C \jjnt_{\Og\times(0,T)} e^{-2s\alpha}
\left[ \lambda^6 (s\phi)^5 |\varphi| + \lambda^5 (s\phi)^4 |\nabla\varphi| + \lambda^4 (s\phi)^3 |\Delta\varphi| \right] |\psi| \;dxdt\\
\disp \phantom{M'_1} \leq \eps I(s,\lambda; \varphi) + C_\eps \jjnt_{\Og\times(0,T)} e^{-2s\alpha} \lambda^8 (s\phi)^7 |\psi|^2 dxdt.
\end{array}
   \end{equation}
   On the other hand,
   \begin{equation}\label{estim-M2p}
M'_2 \leq \eps I(s,\lambda; \varphi) + C_\eps \jjnt_{\Og\times(0,T)} e^{-2s\alpha} \lambda^4 (s\phi)^3 |\psi|^2 dxdt.
   \end{equation}

   From \eqref{5l}, \eqref{new-p} and~\eqref{estim-M1p}--\eqref{estim-M2p}, we find that
   $$
\begin{array}{l}
\disp I(s,\lambda; \varphi) + \tilde{I}(s,\lambda; \psi) \leq C\eps I(s,\lambda; \varphi) + C \jjnt_{\Og\times(0,T)} e^{-2s\alpha} \lambda^4 (s\phi)^3 |\psi|^2 dxdt,
\end{array}
   $$
for all small $\eps>0$.

   This ends the proof.
\Fin

\

   Arguing as in the proof of Corollary~\ref{corol2}, the following can be easily established:

\begin{corollary}\label{corol2p}
   Under the assumptions of Proposition~\ref{t2p}, there exists constants $M'$ and~$K'$, depending on $\Om$, $\Og$, $T$, $\|a\|_{L^\infty(Q)}$, $|b|$, $\|c\|_{L^\infty(Q)}$ and~$\|d\|_{L^\infty}$, such that every solution $(\varphi,\psi)$ to (\ref{as}) verifies:
   \begin{equation}\label{dop}
    \|\varphi(\cdot\,,0)\|^2 + \limsup_{t \to 0^+} \|\psi(\cdot\,,t)\|^2 \leq M' \jjnt_{\Og\times(0,T)}e^{-2K'/(T-t)} |\psi|^2dxdt.
   \end{equation}
\end{corollary}


\section{The null controllability of the linearized systems}\label{Sec3}

   In this Section, we will deduce from the observability estimates \eqref{do} and~\eqref{dop} null controllability results for \eqref{1.1-lin} and~\eqref{1.1p-lin}.

   More precisely, we have:

\begin{theorem}\label{cnl}
   Assume that \eqref{hyp-1.1} holds.
   Then \eqref{1.1-lin} is null-controllable at any time $T > 0$.
   That is to say, for any $y^0 \in L^2(\Om)$ there exist null controls $v \in L^2(\Og\times(0,T))$ for \eqref{1.1-lin} satisfying
   \begin{equation} \label{eq3}
\|v\|_{L^2(\Og\times(0,T))}\leq C \|y^0\|,
   \end{equation}
where $C$ only depends on $\Om$, $\Og$, $T$, $\|a\|_{L^\infty(Q)}$, $\|b\|_{L^\infty}$, $|c|$ and~$\|d\|_{L^\infty}$.
\end{theorem}

\noindent
{\sc Proof:}
   There are several ways to prove that the observability inequality \eqref{do} implies the null controllability of \eqref{1.1-lin}.
   One of them is the following.

   For any $v \in L^2(\Og \times (0,T))$ and any $\eps>0$, let us set
   \begin{equation}\label{fct1}
J_\eps(v) = \jjnt_{\Og\times(0,T)} e^{2K/(T-t)} |v|^2 dxdt+ \frac{1}{\eps} \|y(\cdot\,,T)\|^2 .
   \end{equation}
   Here, $(y, z)$ is the solution to~\eqref{1.1-lin} associated to the initial data $y^0$.
   It is not difficult to check that $v \mapsto  J_\eps(v)$ is lower semi-continuous, strictly convex and coercive in $L^2(\Og \times (0,T))$.
   Hence, it possesses a unique minimizer $v_\eps \in L^2(Q)$.
   We will denote by $(y_\eps, z_\eps)$ the associated state.

   We will show that, at least for a subsequence, $v_\eps$ converges weakly in~$L^2(\Og \times (0,T))$ towards a control $v \in L^2(\Og \times (0,T))$ and the associated $y_\eps(\cdot\,,T)$ converges strongly in~$L^2(\Om)$ to zero.
   Obviously, this proves that $v$ is a null control for \eqref{1.1-lin}, i.e.~that the state associated to $v$ satisfies~\eqref{eq2n}.

   Notice that the unique minimizer of~\eqref{fct1} is characterized by the following optimality system:
   \begin{equation}\label{eq1ep}
\left\{
\begin{array}{l}
\disp y_{\eps,t}-\Delta y_{\eps}= ay_{\eps} + bz_{\eps}+v_\eps1_\Og\ \ \ \mbox{ in } \ Q, \\
\disp -\Delta z_{\eps}=cy_{\eps} + dz_{\eps} \ \ \ \mbox{ in } Q,\\
\disp y_{\eps} = z_{\eps} = 0  \ \ \ \mbox{ on }\ \Sigma,\\
\disp y_{\eps}(x,0) = y^0(x)  \ \ \ \mbox{ in } \ \Omega,
\end{array} \right.
   \end{equation}
   \begin{equation}
\left\{ \begin{array}{l} \label{asep}
\disp -\varphi_{\eps,t}-\Delta \varphi_{\eps}= a\varphi_{\eps} + c\psi_{\eps} \ \ \ \mbox{ in }\ Q  \\
\disp -\Delta \psi_{\eps}= b\varphi_{\eps} + d\psi_{\eps} \ \ \ \mbox{ in } \ Q\\
\disp \varphi_{\eps} = \psi_{\eps}= 0\ \ \ \;\mbox{ on } \ \Sigma \\
\disp \varphi_{\eps}(x,T) =-\frac{1}{\eps}y_{\eps} (x,T)\ \;\;\; \mbox{ in } \;\Omega,
\end{array}\right.
   \end{equation}
   \begin{equation}\label{vep}
v_\eps = \bigl. e^{-2K/(T-t)} \varphi_\eps \bigr|_{\Og \times (0,T)} \,
   \end{equation}
(see for instance~\cite{LT}; see also~\cite{Cara-Guerrero}).

   By multiplying both sides of $(\ref{asep})_1$ by $y_\eps$ and both sides of $(\ref{asep})_2$ by $z_\eps$, integrating in time and space and adding the resulting identities, we obtain:
   \begin{equation}\label{d1}
\begin{array}{l}
\displaystyle \jjnt_{\Og\times(0,T)}  e^{-2K/(T-t)} |\varphi_{\eps}|^2 dxdt + \frac 1\eps \|y_\eps(\cdot\,,T)\|^2 \\[11pt]
\disp \leq \|y^0\|\; \|\varphi_\eps(\cdot\,,0)\|
\leq M^{1/2} \|y^0\| \left( \jjnt_{\Og\times(0,T)} e^{-2K/(T-t)} |\varphi_{\eps}|^2 dxdt \right)^{1/2} \\[11pt]
\disp \leq C\|y^0\|^2 + \frac 12 \jjnt_{\Og\times(0,T)} e^{-2K/(T-t)} |\varphi_{\eps}|^2 dxdt ,
\end{array}
   \end{equation}
where we have used the observability estimate \eqref{do}.

   From (\ref{d1}), we see that
   \begin{equation}\label{d1.1}
\jjnt_{\Og\times(0,T)} e^{2K/(T-t)} |v_{\eps}|^2 dxdt = \jjnt_{\Og\times(0,T)} e^{-2K/(T-t)} |\varphi_{\eps}|^2 dxdt \leq C\|y^0\|^2.
   \end{equation}

   Consequently, at least for a subsequence, one has
   \begin{equation}\label{d3}
\begin{array}{l}
\displaystyle v_{\eps} \to v \ \text{ weakly in } \ L^2(\Og\times (0,T)), \ \text{ with } \ \|v\|_{L^2(\Og\times (0,T))}\leq C \|y^0\| .
\end{array}
   \end{equation}

   We also have from (\ref{d1}) that
   \begin{equation}\label{d4}
\begin{array}{l}
\displaystyle y_{\eps}(\cdot\,,T) \to 0 \ \text{ strongly in }\ L^2(\Om) \ \text{ as } \ \eps \to 0.
\end{array}
   \end{equation}

   From the usual energy method, it is clear that a subsequence can be extracted such that
   \begin{equation}\label{d5}
\begin{array}{l}
\displaystyle y_{\eps} \to y \ \text{ and } \ z_{\eps} \to z \;\;\mbox{strongly in }\;L^2(Q)\;\;\mbox{ as }\;\eps \to 0
\end{array}
   \end{equation}
(see for instance~\cite{9}) and, consequently, the limit $v$ is such that the solution $(y,z)$ to~\eqref{1.1-lin} satisfies \eqref{eq2n}.

   This ends the proof of Theorem~\ref{cnl}.
\Fin

\

   We also deduce from Corollary~\ref{corol2p} the following result:

\begin{theorem}\label{cnlp}
   Assume that \eqref{hyp-1.1p} holds.
   Then \eqref{1.1p-lin} is null-controllable at any time $T > 0$.
   In other words, for any $y^0 \in L^2(\Om)$ there exist null controls $w \in L^2(\Og\times(0,T))$ for \eqref{1.1p-lin} such that
   $$
\|w\|_{L^2(\Og\times(0,T))}\leq C \|y^0\|,
   $$
where $C$ only depends on $\Om$, $\Og$, $T$, $\|a\|_\infty$, $|b|$, $\|c\|_\infty$ and~$\|d\|_\infty$.
\end{theorem}

\noindent
{\sc Proof:}
   It is very similar to the proof of Thoerem~\ref{cnl}.

   Indeed, we can introduce the functional $L_\eps$, with
   \begin{equation}\label{fct1}
L_\eps(v) := \jjnt_{\Og\times(0,T)} e^{2K/(T-t)} |w|^2 \;dxdt+ \frac{1}{\eps} \|y(\cdot\,,T)\|^2
   \end{equation}
and we can again check that $L_\eps$ possesses exactly one minimizer $w_\eps \in L^2(Q)$.
   The optimality system is now
   $$
\left\{
\begin{array}{l}
\disp y_{\eps,t}-\Delta y_{\eps}= ay_{\eps} + bz_{\eps}\ \ \ \mbox{ in } \ Q, \\
\disp -\Delta z_{\eps}=cy_{\eps} + dz_{\eps} + w_\eps1_\Og\ \ \ \mbox{ in } Q,\\
\disp y_{\eps} = z_{\eps} = 0  \ \ \ \mbox{ on }\ \Sigma,\\
\disp y_{\eps}(x,0) = y^0(x)  \ \ \ \mbox{ in } \ \Omega,
\end{array} \right.
   $$
   $$
\left\{ \begin{array}{l}
\disp -\varphi_{\eps,t}-\Delta \varphi_{\eps}= a\varphi_{\eps} + c\psi_{\eps} \ \ \ \mbox{ in }\ Q  \\
\disp -\Delta \psi_{\eps}= b\varphi_{\eps} + d\psi_{\eps} \ \ \ \mbox{ in } \ Q\\
\disp \varphi_{\eps} = \psi_{\eps}= 0\ \ \ \;\mbox{ on } \ \Sigma \\
\disp \varphi_{\eps}(x,T) =-\frac{1}{\eps}y_{\eps} (x,T)\ \;\;\; \mbox{ in } \;\Omega,
\end{array}\right.
   $$
   $$
w_\eps = \bigl. e^{-2K/(T-t)} \psi_\eps \bigr|_{\Og \times (0,T)} \,.
   $$

   We can argue as before and deduce that
   \begin{equation}\label{d1.1}
   \begin{array}{l}
\jjnt_{\Og\times(0,T)} e^{2K/(T-t)} |w_{\eps}|^2 dxdt\\[5pt]\disp
= \jjnt_{\Og\times(0,T)} e^{-2K/(T-t)} |\psi_{\eps}|^2 dxdt \leq C\|y^0\|^2,
\end{array}
   \end{equation}
whence a weakly convergent sequence of control exists and, in the limit, we get a null control for~\eqref{1.1p-lin}.

   Notice that, by construction,
   $$
\limsup_{t \to T^-} \|w(\cdot\,,t)\|_{L^2(\Og)} = 0.
   $$
   This, together with the energy estimates
   $$
\|z_\eps(\cdot\,,t)\|_{H_0^1} \leq C \left(\|y(\cdot\,,t)\| +  \|w(\cdot\,,t)\|_{L^2(\Og)}\right),
   $$
ensures the second part of~\eqref{eq2n}.
\Fin


\section{The null controllability of the semilinear systems}\label{Sec4}

   In this Section, we present the proofs of the main results in this paper, namely Theorems~\ref{t1n} and~\ref{t2n}.
   They will be obtained by combining the linear controllability results in the previous Section and a standard fixed-point argument.

\

\noindent
{\sc Proof of Theorem~\ref{t1n}:}
   Let us first assume that, in~\eqref{1.2a}, $F_0$ is $C^1$.
   In view of~\eqref{1.2a} and~\eqref{1.2b}, we observe that \eqref{1.1} can be written as follows:
   $$
\left\{
\begin{array}{l}
\disp y_{t} - \Delta y = A_0(y) y+bz+ v1_{\Og}\;\;\mbox{in}\;\; Q,\\
\disp -\Delta z=cy+dz\;\;\mbox{in}\;\; Q,\\
\disp y(x,t)=0,\; z(x,t)=0 \;\;\mbox{on}\;\; \Sigma,\\
\disp y(x,0)=y^{0}(x)\;\;\mbox{in}\;\; \Om,
\end{array}
\right.
   $$
with
   \begin{equation}\label{4.1p}
A_0(s) = \left\{\begin{array}{ll}
                \disp {F_0(s) \over s} & \ \text{ if } \ s\not=0, \\
                \noalign{\smallskip}
                \disp F'_0(0)                & \ \text{ otherwise. }
                \end{array}\right.
   \end{equation}

   For any $k \in L^2(Q)$, let us consider the linear system
   \begin{equation}\label{tk1}
\left\{
\begin{array}{l}
 \disp y_{t} - \Delta y =A_0(k) y+bz+ v1_{\Og}\;\;\mbox{in}\;\; Q,\\
\disp -\Delta z=cy+dz\;\;\mbox{in}\;\; Q,\\
\disp
 y(x,t)=0,\; z(x,t)=0 \;\;\mbox{on}\;\; \Sigma,\\
\disp
 y(x,0)=y^{0}(x)\;\;\mbox{in}\;\; \Om.
\end{array}
\right.
   \end{equation}

   In view of Theorem \ref{cnl}, there exists controls $v \in L^2(\Og \times (0,T))$ such that the associated states $(y,z)$ satisfy~\eqref{eq2n} and~\eqref{eq3}, where the constant $C$ only depends on $\Om$, $\Og$, $T$, $\|F_0\|_{C^1(\R)}$, $|b|$, $|c|$ and~$|d|$.

   Let us introduce the mapping $\Phi : L^2(Q) \mapsto 2^{L^2(Q)}$, as follows:
   for any $k \in L^2(Q)$, we set by definition
   $$
\Phi_0(k) = \{\, v \in L^2(\Og \times (0,T)) : \text{ the solution to~\eqref{tk1} satisfies \eqref{eq2n} and~\eqref{eq3}} \,\}
   $$
and
   $$
\Phi(k) = \{\, y \in L^2(Q): \text{ $(y,z)$ solves \eqref{tk1} for some $v \in \Phi_0(k)$} \,\}.
   $$
   Then $\Phi$ satisfies the hypotheses of {\it Kakutani's Fixed-Point Theorem.}

   Indeed, the following holds:

\begin{itemize}

\item For any $k \in L^2(Q)$, $\Phi(k) \subset L^2(Q)$ is a non-empty compact set.
   Furthermore, there exists a fixed compact set $K \subset L^2(Q)$ such that $\Phi(k) \subset K$ for all $k \in L^2(Q)$.

   This is an obvious consequence of \eqref{eq3}, the energy estimates
   $$
\| y_t \|_{L^2(0,T;H^{-1}(\Om))}^2 + \| y \|_{L^2(0,T;H_0^1(\Om))}^2 \leq C \left( \|v\|_{L^2(\Og \times (0,T))}^2 + \|y^0\|^2 \right)
   $$
   (established in the Appendix) and the compactness of the embedding $W \hookrightarrow L^2(Q)$, where
   $$
W = \{\, \xi \in L^2(0,T;H_0^1(\Om)) : \xi_t \in L^2(0,T;H^{-1}(\Om))\,\}.
   $$

\item $\Phi$ is sequentially upper semicontinuous on~$L^2(Q)$, i.e.~if $k_n, k \in L^2(Q)$ and $k_n \to k$ in~$L^2(Q)$, then
   $$
\disp \limsup_{n \to \infty} \left( \sup_{y \in \Phi(k_n)} \jjnt_Q y\,\xi\;dxdt \right) \leq \sup_{y \in \Phi(k_n)} \jjnt_Q y\,\xi \;dxdt \quad \forall \xi \in L^2(Q).
   $$

\end{itemize}

   Therefore, $\Phi$ possesses at least one fixed-point $y$.

   Obviously, $y$ solves, together with some $z$, the semilinear system \eqref{1.1} for some $v \in L^2(\Og \times (0,T))$ and~\eqref{eq2n} holds.

   Let us now assume that the function $F_0$ in~\eqref{1.2a} is only globally Lipschitz-continuous.
   Suppose that
   $$
|F_0(s_1) - F_0(s_2)| \leq L |s_1 - s_2| \quad \forall s_1,s_2 \in \R.
   $$

   Then, we can find $C^1$ functions $F_{01}, F_{02}, \dots$ with the following properties:

\begin{enumerate}

\item $F_{0n} : \R \mapsto \R$ is $C^1$ and globally Lipschitz-continuous., with Lpschitz constant $L$, i.e.
   Suppose that
   $$
|F_{0n}(s_1) - F_{0n}(s_2)| \leq L |s_1 - s_2| \quad \forall s_1,s_2 \in \R, \quad \forall n \geq 1.
   $$

\item $F_{0n} \to F_0$ uniformly on each compact interval $I \subset \R$.

\end{enumerate}

   For each $n \geq 1$, let us consider the system
   \begin{equation}\label{4.10a}
\left\{
\begin{array}{l}
 \disp y_{t} - \Delta y = F_{0n}(y) y+bz+ v1_{\Og}\;\;\mbox{in}\;\; Q,\\
\disp -\Delta z=cy+dz\;\;\mbox{in}\;\; Q,\\
\disp
 y(x,t)=0,\; z(x,t)=0 \;\;\mbox{on}\;\; \Sigma,\\
\disp
 y(x,0)=y^{0}(x)\;\;\mbox{in}\;\; \Om,
\end{array}
\right.
   \end{equation}

   Let $v_n$ be a null control for \eqref{4.10a} satisfying
   $$
\| v_n \|_{L^2(\Og\times(0,T))} \leq C \|y^0\|,
   $$
with $C$ independent of $n$.
   In view of the previous arguments and the properties of $F_{0n}$, such a $v_n$ exists.
   It is clear that, at least for a subsequence, one has
   $$
v_n \to v \ \text{ weakly in } \ L^2(\Og\times(0,T)),
   $$
where $v$ is a null control for \eqref{1.1} again satisfying~\eqref{eq2n}.

   This ends the proof.
\Fin

\

\noindent
{\sc Proof of Theorem~\ref{t2n}:}
   The proof is similar to the previous one, although a little more intrincate.

   Again, let us first assume that, in~\eqref{1.2a} and~\eqref{1.2c}, the functions $F_0$ and~$f_0$ are $C^1$.
   Then, \eqref{1.2} can be written in the form
   $$
\left\{
\begin{array}{l}
\disp y_{t} - \Delta y = A_0(y) y+bz \;\;\mbox{in}\;\; Q,\\
\disp -\Delta z=C_0(y)y + dz + w1_{\Og}\;\;\mbox{in}\;\; Q,\\
\disp y(x,t)=0,\; z(x,t)=0 \;\;\mbox{on}\;\; \Sigma,\\
\disp y(x,0)=y^{0}(x)\;\;\mbox{in}\;\; \Om,
\end{array}
\right.
   $$
where $A_0$ is given by~\eqref{4.1p} and
   $$
C_0(s) = \left\{\begin{array}{ll}
                \disp {f_0(s) \over s} & \ \text{ if } \ s\not=0, \\
                \noalign{\smallskip}
                \disp f'_0(0)                & \ \text{ otherwise. }
                \end{array}\right.
   $$

   For each $k \in L^2(Q)$, we can consider the linear system
   \begin{equation}\label{4.10c}
\left\{
\begin{array}{l}
\disp y_{t} - \Delta y =A_0(k) y+bz \;\;\mbox{in}\;\; Q,\\
\disp -\Delta z=C_0(k)y + dz + w1_{\Og}\;\;\mbox{in}\;\; Q,\\
\disp y(x,t)=0,\; z(x,t)=0 \;\;\mbox{on}\;\; \Sigma,\\
\disp y(x,0)=y^{0}(x)\;\;\mbox{in}\;\; \Om.
\end{array}
\right.
   \end{equation}

   As in the proof of Theorem~\ref{t1n}, we can prove that there exist controls $w \in L^2(\Og \times (0,T))$ such that the associated solutions to~\eqref{4.10c} satisfy~\eqref{eq2n} and
   $$
\| w \|_{L^2(\Og\times(0,T))} \leq C \|y^0\|,
   $$
for some fixed $C$.
   We can again introduce a multi-valued mapping $\Psi : L^2(Q) \mapsto 2^{L^2(Q)}$
   (similar to~$\Phi$) and we can show that the assumptions of Kakutani's Theorem are satisfied by~$\Psi$, etc.

   An easy adaptation of the remaining results leads to the desired controllability result.
   We omit the details, that can be checked easily.
\Fin


\section{An asymptotic controllability property}\label{Sec5}

   In this Section, we prove that, under appropriate conditions on~$F$ and~$f$
   (in fact the same in~Theorem~\ref{t2n}), the semilinear parabolic system
   \begin{equation}\label{1epsn}
\left\{
\begin{array}{l}
 \disp y_{t} - \Delta y =F(y,z)\;\;\mbox{in}\;\; Q,\\
\disp \eps\,z_t -\Delta z=f(y,z)+ w1_{\Og}\;\;\mbox{in}\;\; Q,\\
\disp
 y(x,t)=0,\; z(x,t)=0 \;\;\mbox{on}\;\; \Sigma,\\
\disp
 y(x,0)=y^{0}(x), \,\,z(x,0)=z^0(x)\;\;\mbox{in}\;\; \Omega,
\end{array}
\right.
   \end{equation}
is uniformly null-controllable as $\eps \to 0$.
   We also prove the convergence of the null controls to a null control for the similar parabolic-elliptic system
   \begin{equation}\label{3epsn*}
\left\{
\begin{array}{l}
 \disp y_{t} - \Delta y =F(y,z)\;\;\mbox{in}\;\; Q,\\
\disp -\Delta z=f(y,z)+ w1_{\Og}\;\;\mbox{in}\;\; Q,\\
\disp y(x,t)=0,\; z(x,t)=0 \;\;\mbox{on}\;\; \Sigma,\\
\disp y(x,0)=y^{0}(x)\;\;\mbox{in}\;\; \Omega.
\end{array}
\right.
   \end{equation}

   To this end, we will first consider the linear system
   \begin{equation}\label{1epsl}
\left\{
\begin{array}{l}
 \disp y_{t} - \Delta y =ay+bz\;\;\mbox{in}\;\; Q,\\
\disp \eps z_t-\Delta z=cy+dz+ w1_{\Og}\;\;\mbox{in}\;\; Q,\\
\disp y(x,t)=0,\; z(x,t)=0 \;\;\mbox{on}\;\; \Sigma,\\
\disp y(x,0)=y^{0}(x), \;\;\;z(x,0)=z^0(x)\;\;\mbox{in}\;\; \Om ,
\end{array}
\right.
   \end{equation}
and we will establish a uniform null controllability result.

   More precisely, the following holds:

\begin{theorem}\label{cnepsl}
   Assume that \eqref{hyp-1.1p} holds.
   Then, for any~$\eps >0$ and any~$y^0, z^0 \in L^2(\Om)$, there exist controls $w_\eps \in L^2(\Og\times(0,T))$ such that the corresponding solutions $(y_\eps,z_\eps)$ to~\eqref{1epsl} satisfy
   $$
y_\eps(x,T)=0,\;\;z_\eps(x,T)=0 \;\;\;\mbox{in}\;\;\Om,
   $$
with an estimate of the form
   \begin{equation}
\|w_\eps\|_{L^2(\Og\times(0,T))}\leq C \left(\|y^0\|+\eps \; \|z^0\| \right). \label{76epsl}
   \end{equation}
where $C$ is independent of~$\eps$.
\end{theorem}

\noindent
{\sc Sketch of the proof:}
   Let us consider the adjoint system of~(\ref{1epsl}), that is:
   $$
\left\{
\begin{array}{l}
- \varphi_t-\Delta \varphi=a\varphi+c\psi\;\;\mbox{in}\;\; Q,\\[5pt]
\disp-\eps \psi_t -\Delta \psi=b\varphi+d\psi\;\;\mbox{in}\;\; Q,\\[5pt]
\disp \varphi=0, \; \psi=0 \;\;\mbox{on}\;\; \Sigma,\\[5pt]
\disp \varphi(x,T)=\varphi^0(x), \,\,\,\psi(x,T)=\psi^0(x)\;\;\mbox{in}\;\; \Om,
\end{array}
\right.
   $$
with $\varphi^0, \psi^0 \in L^2(\Om)$.

   As in Section~\ref{Sec2}, we will use an abridged notation for the weighted integrals concerning $\varphi$, $\psi$ and their derivatives:
   for any positive $\lambda$ and $s$, we set
   $$
\begin{array}{l}
\disp I(s,\lambda, \varphi)=\displaystyle \jjnt_Q \!\!e^{-2s\alpha}\!\left[(s\phi)^{-1}\!\left(|\varphi_t|^2\!+\!|\Delta \varphi|^2 \right)\!+\! \lambda^2 s \phi|\nabla \varphi|^2 \!+\! \lambda^4(s\phi)^3|\varphi|^2\right]\! \,dx\,dt,\\[9pt]
\disp I_\eps(s,\lambda, \psi)=\displaystyle \eps^2 \jjnt_Q e^{-2s\alpha}(s\phi)^{-1}|\psi_t|^2\,dx\,dt\\[9pt]
\disp \phantom{I_\eps(s,\lambda, \psi)} + \jjnt_Q \!\!e^{-2s\alpha}\!\left[(s\phi)^{-1} |\Delta \varphi|^2 \!+\! \lambda^2 s \phi|\nabla \psi|^2 \!+\! \lambda^4(s\phi)^3|\varphi|^2\right]\! \,dx\,dt
\end{array}
   $$
   Then one has
   \begin{equation}\label{60epsl}
\disp I(s,\lambda,\varphi)+I_\eps(s,\lambda, \psi)\leq C \jjnt_{\Og\times(0,T)} (s\phi)^7 \lambda^8e^{-2s\al(x,t)}|\psi|^2 \,dx\,dt,
   \end{equation}
where $C$ is independent of $\eps$.

   The proof of \eqref{60epsl} is very similar to the proof of Proposition~\ref{t2} and will be omitted.
   An almost immediate consequence is the following observability inequality:
   \begin{equation}\label{DOeps}
\|\varphi(\cdot\,,0)\|^2+\eps \;\|\psi(\cdot\,,0)\|^2 \leq C \jjnt_{\Og\times(0,T)} e^{2s\alpha}(s\phi)^7 \ld^8|\psi|^2\,dx\,dt.
   \end{equation}

   Now, arguing as in~Section~\ref{Sec3}, it becomes clear that the uniform null controllability property of (\ref{1epsl}) is implied by the observability estimate (\ref{DOeps}).
\Fin

\

   From this result, we get the following for the semilinear systems \eqref{1epsn} and~\eqref{3epsn*}:

\begin{theorem}\label{TFep}
   Under the assumptions of Theorem~\ref{t2n}, for any~$\eps >0$ and any~$y^0, z^0 \in L^2(\Omega)$, there exist null controls $w_\eps \in L^2(\Og \times(0,T))$ for~\eqref{1epsn}.
   They can be chosen such that, at least for a subsequence,
   $$
w_\eps \to w \;\;\;\mbox{ weakly in }\;\;L^2(\Og \times (0,T)),
   $$
where $w$ is a null control for \eqref{3epsn*}.
\end{theorem}

\noindent
{\sc Sketch of the proof:}
   For the proof of the first assertion, it suffices to argue as in the proof of Theorem~\ref{t2n}.

   Indeed, for any $\eps > 0$ and any fixed $k \in L^2(Q)$, we can apply Theorem~\ref{cnepsl} to the linearized system
   \begin{equation}\label{6epsn}
\left\{
\begin{array}{l}
 \disp y_{t} - \Delta y =A_0(k)y+b z\;\;\mbox{in}\;\; Q,\\
\disp \eps\,z_t-\Delta z=C_0(k)y+dz+ w1_{\Og}\;\;\mbox{in}\;\; Q,\\
\disp y(x,t)=0,\; z(x,t)=0 \;\;\mbox{on}\;\; \Sigma,\\
\disp  y(x,0)=y^{0}(x), \,\,z(x,0)=z^0(x)\;\;\mbox{in}\;\; \Omega.
\end{array}
\right.
   \end{equation}

   We deduce that, for any $\eps >0$, (\ref{6epsn}) is null-controllable, with controls $w_\eps$ satisfying (\ref{76epsl})
   (where $C$ is independent of~$\eps$).
   Observe that, in fact, we can get a stronger estimate:
   $$
\jjnt_{\Og\times(0,T)} e^{2K/(T-t)} |w_{\eps}|^2 dxdt\leq C\|y^0\|^2.
   $$

   We can again introduce a multi-valued mapping $\Psi_\eps$
   (similar to~$\Psi$) and we can show that the assumptions of Kakutani's Theorem are satisfied by~$\Psi_\eps$.
   Therefore, there exist controls $w_\eps \in L^2(\Og\times(0,T))$ satisfying (\ref{76epsl}) such that the associated solutions to~(\ref{1epsn}) satisfy
   $$
y_\eps(x,T)=0,\;\;z_\eps(x,T)=0 \;\;\;\mbox{in}\;\;\Om.
   $$

   Let us multiply both sides of $(\ref{6epsn})_1$
   (resp.~$(\ref{6epsn})_2$) by $\disp y_\eps$ (resp.~$\disp z_\eps$) and let us integrate in $Q$.
   We easily obtain the following for all $t$:
   $$
\begin{array}{l}
\disp \|y_\eps(\cdot\,,t)\|^2 + \eps \|z_\eps(\cdot\,,t)\|^2 + \jnt_0^t \left( \|y_\eps(\cdot\,,s)\|_{H_0^1}^2 + \|y_\eps(\cdot\,,s)\|_{H_0^1}^2 \right) \,ds \\
\disp \ \ \leq \|y^0\|^2 + \eps \|z^0\|^2 + C \jjnt_{\Og\times(0,t)} |w_\eps|^2 \,dx\,ds + C \jnt_0^t \|y_\eps(\cdot\,,s)\|^2 \,ds.
\end{array}
   $$

   Using Gronwall's inequality and extracting appropriate subsequences, we deduce that, at least for a subsequence, $w_\eps$, $y_\eps$ and $z_\eps$ respectively converge to $w$, $y$ and $z$, where $(y,z)$ solves (\ref{3epsn*}) and satisfies (\ref{eq2n}).
\Fin


\section{Appendix A: Some technical results}\label{Sec6}

   In this Appendix, for completeness, we give a theoretical result for the semilinear system
   \begin{equation}\label{a1}
\left\{
\begin{array}{l}
\disp y_{t} - \Delta y = F(y,z)+ v1_{\Og}\;\;\mbox{in}\;\; Q,\\
\disp -\Delta z= f(y,z)+ w1_{\Og}\;\;\mbox{in}\;\; Q,\\
\disp y(x,t)=0,\; z(x,t)=0 \;\;\mbox{on}\;\; \Sigma,\\
\disp y(x,0)=y^{0}(x)\;\;\mbox{in}\;\; \Om,
\end{array}
\right.
   \end{equation}
where $F$ and $f$ satisfy (\ref{1.2a}) and (\ref{1.2c}), respectively.

   We have the following:

\begin{theorem} \label{tha1}
   For any given $y^0 \in L^2(\Om)$ and $v, w \in L^2(\Og \times (0,T))$, \eqref{a1} possesses a unique solution $(y,z)$, with
   \begin{equation}\label{rea1}
y \in L^2(0,T;H_0^1(\Om)), \quad y_t \in L^2(0,T;H^{-1}(\Om)), \quad z \in L^2(0,T;D(-\Delta)).
   \end{equation}
\end{theorem}


\noindent
{\sc Sketch of the proof:}
   We will apply the Faedo-Galerkin method,  see~J.-L.~Lions~\cite{9}.

   Let $\{h_j\}$ be a {\it special} basis of~$H_0^1(\Om)$, more precisely, the basis formed by the eigenfuncions of the Dirichlet Laplacian in~$\Om$.
   Let us introduce the finite-dimensional Galerkin approximations as follows:
   find $y_N, z_N$, with $y_N(t), z_N(t) \in V_N$ for all $t$, such that
   \begin{eqnarray}&
(y_N'(t),h) +(\nabla y_N(t),\nabla h) = (F(y_N,z_N),h)+(v1_\Og,h) \quad \forall h \in V_N, \label{a2}
\\[4pt]&
\hspace{-10pt}(\nabla z_N(t),\nabla k) = (f(y_N,z_N),k)+(w1_\Og,k) \quad \forall k \in V_N, \label{a3}
\\[4pt]&
\begin{array}{l}
y_N(0) = y_{0N}.
\end{array} \label{a4}
   \end{eqnarray}

   Here, $V_N = [h_1,h_2,\dots,h_N]$ is the subspace of $H_0^1(\Om)$ spanned by the first~$N$ eigenfunctions $h_j$ and $y_{0N} = P_N y_0$,  where $P_N: L^2(\Om) \mapsto V_N$ is the orthogonal projector.

   This Cauchy problem has at least one local solution on $[0,t_N)$.
   That the maximal solution is defined in the whole interval $[0,T]$ is a consequence of the estimates given below.

\

{\sc Estimates~I :}
   Let us set $h = y_N(t)$ in~(\ref{a2}) and $k=z_N(t)$ in~(\ref{a3}).
   After some computations we obtain the following for any~$t$ and for all small $\delta > 0$:
   \begin{equation}\label{a5}
\begin{array}{l}
\disp {1\over2} \|y_N(t)\|^2 + \jnt_0^t \|y_N(s)\|_{H_0^1}^2 \,ds + \left( 1 - {\mu + \delta \over \mu_1} \right) \jnt_0^t \|z_N(s)\|_{H_0^1}^2 \,ds \\
\disp \ \ \leq {1\over2} \|y^0\|^2 + C_\delta \jjnt_{\Og\times(0,t)} \left( |v|^2 + |w|^2 \right) \,dx\,ds + C_\delta \jnt_0^t \|y_N(s)\|^2 \,ds
\end{array}
   \end{equation}
and
   \begin{equation}\label{a5-z}
\begin{array}{l}
\disp\left( 1 - {\mu + \delta \over \mu_1} \right) \|z_N(t)\|_{H_0^1} \leq C_\delta \left( \|y_N(t)\| + \|w\|_{L^2(\Og)} \right).
\end{array}
   \end{equation}


\

{\sc Estimates~II }
   Since the $h_j$ are the eigenfunctions of $-\Delta$ in~$H_0^1(\Om)$, one has
   \begin{equation}
\|y'_{N}(t)\|_{H^{-1}(\Omega)} \leq \| (\Delta y_N + F(y_N,z_N) + v1_\Og)(t) \|_{H^{-1}(\Omega)}
   \end{equation}
for all $t$.
   Therefore, we get the following for some $C>0$:
    \begin{equation}\label{a6}
\|y'_{N}(t)\|_{H^{-1}(\Omega)} \leq C \left( \|y_N(t)\|_{H_0^1} + \|(v1_\Og)(\cdot\,,t)\| + \|(w1_\Og)(\cdot\,,t)\| \right).
   \end{equation}

\

   From \eqref{a5}--\eqref{a6}, it is standard to deduce that, at least for a subsequence, the $y_N$ and~$z_N$ converge to a solution to~\eqref{a1}.
   This solution must satisfy
   $$
y \in L^2(0,T;H_0^1(\Om)), \quad y_t \in L^2(0,T;H^{-1}(\Om)), \quad z \in L^2(0,T;H_0^1(\Om)).
   $$

   Furthermore, from the usual elliptic estimates, we also obtain that $z \in L^2(0,T;D(\Delta))$.
   This yields~\eqref{rea1}.

   The uniqueness of the solution is also a standard consequence of the previous estimates
   (written for $y := y^1 - y^2$ and $z := z^1 - z^2$ where $(y^1,z^1)$ and~$(y^2,z^2)$ are assumed to solve the system) and the global Lipschitz-continuity of $F$ and~$f$.
\Fin


{\footnotesize

\end{document}